\RequirePackage{fix-cm}
\documentclass[smallextended]{svjour3}       % onecolumn (second format)
\smartqed  % flush right qed marks, e.g. at end of proof
\usepackage{amsmath, amscd, amsfonts, amssymb, graphicx, color}
\begin{document}

\title{Operator inequalities associated with the Kantorovich type inequalities for $s$-convex functions %\thanks{Grants or other notes
}

\titlerunning{Operator inequalities associated with the Kantorovich type inequalities}        % if too long for running head

\author{Ismail Nikoufar and Davuod Saeedi}


\institute{I. Nikoufar (Corresponding Author) \at
              Department of Mathematics, Payame Noor University, P.O. Box 19395-4697 Tehran, Iran \\
            %  Tel.: +123-45-678910\\
             % Fax: +123-45-678910\\
              \email{nikoufar@pnu.ac.ir}           %  \\
\\
D. Saeedi \at
              Department of Mathematics, Payame Noor University, P.O. Box 19395-4697 Tehran, Iran  \\
            %  Tel.: +123-45-678910\\
             % Fax: +123-45-678910\\
              \email{dsaeedi3961@gmail.com}           %  \\
}

\date{Received: date / Accepted: date}

\maketitle

\begin{abstract}
In this paper,
we prove some operator inequalities associated
with an extension of the Kantorovich type inequality for $s$-convex functions.
We also give an application to the order preserving power inequality of three variables
and find a better lower bound for the numerical radius of a positive operator
under some conditions.

\end{abstract}

\section{Introduction}

The class of $s$-convex functions in the second sense is defined in \cite{Hudzik-AEqMath1994}.
Note that all $s$-convex functions in the second sense are nonnegative \cite{Hudzik-AEqMath1994}.
There is an identity between the class of 1-convex functions and the class of ordinary convex functions.
Indeed, the $s$-convexity means just the convexity when $s = 1$, no mater in the first sense or in the second sense.
Moreover, when $s\to 0$ we reach another class of functions the so-called $P$-class functions.
We refer to see some inequalities for $P$-class functions in \cite{Nikoufar-Filomat2020}.
Some properties of $s$-convex functions in both senses are considered in \cite{Hudzik-AEqMath1994} where various examples and counterexamples are given.
For many research results and applications related to the $s$-convex functions see
\cite{Dragomir-DM1999,Kirmaci-AMC2007,Pinheiro-AEqMath2007,Set-CA2012,Dragomir-book2002,Nikoufar-CAOT2021} and so on.

A function $f : [0,\infty) \to\Bbb{R}$ is $s$-convex in the first sense if
\begin{equation}\label{s-convex}
f(\alpha x + \beta y)\leq \alpha^sf(x) + \beta^sf(y)
\end{equation}
holds for all $x, y \in [0,\infty)$ and $\alpha,\beta\geq 0$ with $\alpha^s+\beta^s=1$ and for some fixed $s \in (0, 1]$.
A function $f : [0,\infty) \to\Bbb{R}$ is $s$-convex in the second sense whenever
\begin{equation}\label{s-convex}
f(\lambda x + (1-\lambda)y)\leq \lambda^sf(x) + (1-\lambda)^sf(y)
\end{equation}
holds for all $x, y \in [0,\infty)$ and $\lambda\in [0, 1]$ and for some fixed $s \in (0, 1]$.

% %%%%%%%%%%%%%%%%%%%%%%%%%%%%%%%%%%%%%%%%%%%%%%%%%%%%%%%%
%
\begin{lemma}\cite{Nikoufar-CAOT2021}\label{lm-scnx}
Let $0<s\leq r<1$. Then, the function $f(t)=t^{r}$, $t\geq 0$ is $s$-convex in the second sense.
\end{lemma}
\begin{lemma}\cite{Nikoufar-CAOT2021}\label{lm-scnx-norm}
Let $X$ be a normed space and $0<s\leq r<1$. Then, the function $f:X\to\Bbb{R}$ defined by $f(t)=||t||^{r}$ is $s$-convex in the second sense.
\end{lemma}

Let $\mathcal{H}$ be a Hilbert space and $B(\mathcal{H})$ the algebra of all bounded linear operators on
$\mathcal{H}$. An operator $A$ in $B(\mathcal{H})$ is positive if $\langle Ax, x\rangle \geq 0$
for all $x\in\mathcal{H}$ and denoted by $A \geq 0$.
By $Sp(A)$ we denote the spectrum of an operator $A\in B(\mathcal{H})$.

The celebrated Kantorovich inequality asserts that
if $A$ is an operator on $\mathcal H$ such
that $M\geq A\geq m>0$, then
\begin{equation}\label{KI}
\langle A^{-1}x,x\rangle\langle Ax,x\rangle\leq \frac{(m+M)^2}{4mM}
\end{equation}
holds for every unit vector $x$ in $\mathcal H$.
We call the constant $\frac{(m+M)^2}{4mM}$ Kantorovich constant. The inequality \eqref{KI}
is closely related to properties of convex
functions, and many authors have given many results and comments \cite{Furuta-TF2001,Furuta-Zagreb2005,Mond-HJM1993,Mond-IJM1993,Takahashi-MJ1999}.
By replacing $x$ with $\frac{A^{\frac{1}{2}}x}{||A^{\frac{1}{2}}x||}$ in \eqref{KI}
it is well known
that the inequality \eqref{KI} is equivalent to the following inequality
$$
\langle A^2x,x\rangle \leq \frac{(m+M)^2}{4mM} \langle Ax,x\rangle^2
$$
for every unit vector $x$ in $\mathcal H$.

We proved in \cite{Nikoufar-CAOT2021} some inequalities for
self-adjoint operators on a Hilbert space including
an operator version of the Jensen's inequality
for $s$-convex functions as follows.

\begin{theorem} \label{CAOT-thm1}
Let $A$ be a self-adjoint operator on
the Hilbert space $\mathcal{H}$ and assume that $Sp(A)\subseteq [m,M]$ for some scalars $m,M$ with
$0<m < M$. If $f$ is a continuous $s$-convex function on $[m, M]$, then
\begin{equation}\label{eq0-thm1}
f(\langle Ax, x\rangle) \leq 2^{1-s}\langle f(A)x,x\rangle
\end{equation}
for each $x\in \mathcal{H}$ with $||x||=1$. %$\langle x,x\rangle = 1$.
\end{theorem}

%As a consequence of this result, we proved the following result.

%\begin{theorem} \label{CAOT-thm5}
%\cite[Theorem 5]{Nikoufar-CAOT2021}
%Let $A$ be a self-adjoint operator on
%the Hilbert space $\mathcal{H}$ and assume that $Sp(A)\subseteq [m,M]$ for some scalars $m,M$ with
%$0<m < M$. If $f$ is a continuous $s$-convex function on $[m, M]$, then
%\begin{equation}\label{eq0-thm3}
%\langle f(A)x, x\rangle\leq 2^{1-s}\Big(\frac{M-\langle Ax,x\rangle}{M-m}f(m)+\frac{\langle Ax,x\rangle-m}{M-m}f(M)\Big)
%\end{equation}
%for each $x\in \mathcal{H}$ with $||x||=1$. %$\langle x,x\rangle = 1$.
%\end{theorem}

Furuta \cite{Furuta-JIA1998} gave some operator inequalities associated
with an extension of the Kantorovich inequality for convex function.
%We motivated to prove some operator inequalities associated
%with an extension of the Kantorovich type inequality for $s$-convex functions.
%Throughout this paper, we use $s$-convex functions in the second sense.
In Section 2, we prove some operator inequalities associated
with an extension of the Kantorovich type inequality for $s$-convex functions.
We also give an application to the order preserving power inequality of three variables
and provide a better lower bound for the numerical radius of a positive operator
under some conditions.
In Section 3, we show a difference version of Furuta and Giga's result for $s$-convex functions
and provide some applications.

We remark that in the special case when $s=1$, our results are valid for the convex function.
%%%% However, it is essential to note that
Throughout this paper, we consider $s$-convex functions in the second sense for some fixed $s \in (0, 1]$.

\section{Operator inequalities associated with the Kantorovich type inequality}
%In this section, we give some operator inequalities associated
%with an extension of the Kantorovich type inequality for $s$-convex functions
%and we provide an application to the order preserving power inequality of three variables.

The following lemma was proved in \cite{Furuta-JIA1998} for the case $q>1$ or $q<0$.

\begin{lemma}\label{lemma1}
\cite[Lemma 1.3]{Furuta-JIA1998}
Let $h(t)$ be defined by
\begin{equation}\label{eq-pp0}
h(t)=\frac{1}{t^q}\bigg(k+\frac{K-k}{M-m}(t-m)\bigg)
\end{equation}
on $[m,M]$, $(M>m>0)$,
and any real numbers $q$, $K$, and $k$. Then
the function $h(t)$ has the following upper bound
on $[m,M]$:
\begin{equation}\label{eq-pp}
\frac{mK-Mk}{(q-1)(M-m)}\bigg(\frac{(q-1)(K-k)}{q(mK-Mk)}\bigg)^q,
\end{equation}
where $m, M, k, K$ and $q$
satisfy any one (i) and (ii) respectively:
\begin{itemize}
\item[(i)] $K>k$, $\frac{K}{M}>\frac{k}{m}$, and $\frac{k}{m}q\leq\frac{K-k}{M-m}\leq\frac{K}{M}q$ holds
for a real number $q>1$,
\item[(ii)] $K<k$, $\frac{K}{M}<\frac{k}{m}$, and $\frac{k}{m}q\leq\frac{K-k}{M-m}\leq\frac{K}{M}q$ holds
for a real number $q<0$.
\end{itemize}
\end{lemma}

%%%%%%%%%%%%%%%
We extend this lemma for the case $0<q<1$ as follows.

\begin{lemma}\label{lemma2}
Let $h(t)$ be defined by \eqref{eq-pp0}
on $[m,M]$, $(M>m>0)$. Then
the function $h(t)$ has the lower bound \eqref{eq-pp}
on $[m,M]$, where $m, M, k, K$ and $q$
satisfy the following conditions:
\begin{itemize}
\item[] $K>k$, $\frac{K}{M}<\frac{k}{m}$, and $\frac{K}{M}q\leq\frac{K-k}{M-m}\leq\frac{k}{m}q$
holds for a real number $0<q<1$.
\end{itemize}
Moreover, in this case the upper bound occurs at the beginning or end of the interval $[m,M]$
and the upper bound value is $\max\{\frac{k}{m^q}, \frac{K}{M^q}\}$.
\end{lemma}
\begin{proof}
By an easy differential calculus, $h'(t_1)=0$
when
$$
t_1=\frac{q}{q-1}\frac{mK-Mk}{K-k}.
$$
On the other hand, $h''(t_1)=\frac{q(Mk-mK)}{(M-m)t_1^{q+2}}>0$ and so $t_1$ gives
the lower bound \eqref{eq-pp} of $h(t)$ on $[m,M]$.
Moreover, $t_1$ is the global minimum point of $h$
and $m$ and $M$ are the local maximum points of $h$ with the maximum values $\frac{k}{m^q}$ and $\frac{K}{M^q}$, respectively.
\end{proof}

We now give some operator inequalities associated with an extension of the Kantorovich inequality.
One target of the current paper is to extend the
domain of $q$ to values inside the interval $(0, 1)$.
We deal with it in Theorem \ref{main-1}.

\begin{theorem}\label{main-0}
Let $A$ be a positive operator on a Hilbert space
$\mathcal H$ satisfying $M\geq A\geq m>0$.
Also let $f(t)$ be a real valued continuous $s$-convex function on $[m,M]$.
Then the following inequality holds for every unit vector $x$ and for any
real number $q$ depending on (i) or (ii):
$$
\langle f(A)x,x\rangle\leq 2^{1-s} K_f(m,M,q)\langle Ax,x\rangle^q,
$$
where
$$
K_f(m,M,q)=\frac{mf(M)-Mf(m)}{(q-1)(M-m)}\bigg(\frac{(q-1)(f(M)-f(m))}{q(mf(M)-Mf(m))}\bigg)^q
$$
under any one of the following conditions (i) and (ii) respectively:
\begin{enumerate}
\item[(i)] $f(M)>f(m)$, $\frac{f(M)}{M}>\frac{f(m)}{m}$ and $\frac{f(m)}{m}q\leq\frac{f(M)-f(m)}{M-m}\leq\frac{f(M)}{M}q$
holds for a real number $q>1$,
\item[(ii)] $f(M)<f(m)$, $\frac{f(M)}{M}<\frac{f(m)}{m}$ and $\frac{f(m)}{m}q\leq\frac{f(M)-f(m)}{M-m}\leq\frac{f(M)}{M}q$
holds for a real number $q<0$.
\end{enumerate}
\end{theorem}
\begin{proof}
Consider
$D=\Big(
\begin{array}{ll}
m & 0\\
0 & M
\end{array}
 \Big)$
 and
 $x=\Bigg(
\begin{array}{l}
\sqrt{\frac{M-t}{M-m}}\\
\sqrt{\frac{t-m}{M-m}}
\end{array}
 \Bigg)$.
According to Theorem \ref{CAOT-thm1}, one has
\begin{align}\label{eq-qq1}
f(t)&=f(\langle Dx,x \rangle)\nonumber\\
&\leq 2^{1-s}\langle f(D)x,x\rangle\nonumber\\
&=2^{1-s}\Big(\frac{M-t}{M-m}f(m)+\frac{t-m}{M-m}f(M)\Big)\nonumber\\
&=2^{1-s}\Big(f(m)+\frac{f(M)-f(m)}{M-m}(t-m)\Big)
\end{align}
for every $t\in [m,M]$.
So,
\begin{align}\label{eq-qq2}
\langle f(A)x,x\rangle\leq 2^{1-s}\Big(f(m)+\frac{f(M)-f(m)}{M-m}(\langle Ax,x\rangle-m)\Big)
\end{align}
for every unit vector $x$.
Multiplying $\langle Ax,x\rangle^{-q}$ on both sides of \eqref{eq-qq2}, one can deduce
\begin{align}\label{eq-qq3}
\langle Ax,x\rangle^{-q}\langle f(A)x,x\rangle
&\leq 2^{1-s} h(\langle Ax,x\rangle),
\end{align}
where
\begin{align}\label{eq-qq4}
h(t)=t^{-q}\Big(f(m)+\frac{f(M)-f(m)}{M-m}(t-m)\Big)
\end{align}
for every $t\in [m,M]$.
Consequently, it follows from \eqref{eq-qq3} that
\begin{align}\label{eq-qq5}
\langle f(A)x,x\rangle\leq 2^{1-s}\Big(\max_{m\leq t\leq M}h(t)\Big)\langle Ax,x\rangle^{q}.
\end{align}
Putting $K=f(M)$ and $k=f(m)$ in Lemma \ref{lemma1}, then the conditions in Theorem \ref{main-0}
corresponds to the conditions in Lemma \ref{lemma1}, so the proof is complete by \eqref{eq-qq5}.
\end{proof}
In the special case when we consider $s=1$, Theorem \ref{main-0} is valid for the convex functions
and recovers \cite[Theorem 1.1]{Furuta-JIA1998}.

\begin{lemma}\label{lemma-log}
The function $f(t)=\log(t)$ on $[m,M]$ with $m\geq 1$ is $s$-convex. % in the second sense for some $0<s<1$.
\end{lemma}
\begin{proof}
Let $x,y\in [m,M]$, $1<x<y$ and $0<\alpha< 1$.
We realize two cases:

(I) If $\alpha\leq 1-\alpha$, then by weighted AM-GM inequality \cite{Cvetkovski} one has
$$
x^{\alpha}y^{\alpha}\leq x^{\alpha}y^{1-\alpha}\leq \alpha x+(1-\alpha)y.
$$
On the other hand, one has
$$
\alpha x+(1-\alpha)y\leq y< xy.
$$
By combining these inequalities one can get
$$
x^{\alpha}y^{\alpha}\leq  \alpha x+(1-\alpha)y\leq y< xy.
$$
We have
$$
\alpha\log(xy)\leq \log(\alpha x+(1-\alpha)y)\leq \log(y)< \log(xy)
$$
so that
$$
\alpha\leq \frac{\log(\alpha x+(1-\alpha)y)}{\log(xy)}\leq \frac{\log(y)}{\log(xy)}< 1.
$$
It follows that
$$
\log(\alpha)\leq \log\bigg(\frac{\log(\alpha x+(1-\alpha)y)}{\log(xy)}\bigg)\leq \log\bigg(\frac{\log(y)}{\log(xy)}\bigg)< 0,
$$
which entails
$$
1\geq \frac{\log\bigg(\frac{\log(\alpha x+(1-\alpha)y)}{\log(xy)}\bigg)}{\log(\alpha)}\geq \frac{\log\bigg(\frac{\log(y)}{\log(xy)}\bigg)}{\log(\alpha)}> 0.
$$
Define
\begin{align*}
%\theta_1&=\max_{x,y\in [m,M]}\frac{\log\bigg(\frac{\log(y)}{\log(xy)}\bigg)}{\log(\alpha)},\\
\theta&=\min_{x,y\in [m,M], \alpha\in [0,\frac{1}{2}]}\frac{\log\bigg(\frac{\log(\alpha x+(1-\alpha)y)}{\log(xy)}\bigg)}{\log(\alpha)}.
\end{align*}
Choose $s>0$ such that $s\leq \theta$.
This ensures that
$$
s\leq \frac{\log\bigg(\frac{\log(\alpha x+(1-\alpha)y)}{\log(xy)}\bigg)}{\log(\alpha)}
$$
and so
$$
\alpha^s\geq \frac{\log(\alpha x+(1-\alpha)y)}{\log(xy)}.
$$
Hence,
$$
\log(x^{\alpha^s}y^{\alpha^s})\geq \log(\alpha x+(1-\alpha)y).
$$
From (I), we deduce
$$
\log(x^{\alpha^s}y^{(1-\alpha)^s})\geq\log(x^{\alpha^s}y^{\alpha^s}),
$$
which implies
$$
\alpha^s\log(x)+(1-\alpha)^s\log(y)=\log(x^{\alpha^s}y^{(1-\alpha)^s})\geq \log(\alpha x+(1-\alpha)y),
$$

(II) If $\alpha\geq 1-\alpha$, then by weighted AM-GM inequality we get
$$
x^{1-\alpha}y^{1-\alpha}\leq x^{\alpha}y^{1-\alpha}\leq \alpha x+(1-\alpha)y.
$$
We follow a similar path to the proof of the case (I). So,
$$
(1-\alpha)\log(xy)\leq \log(\alpha x+(1-\alpha)y)\leq \log(y)< \log(xy).
$$
Define
%$$
%\theta_1=\max_{x,y\in [m,M]}\frac{\log\bigg(\frac{\log(y)}{\log(xy)}\bigg)}{\log(1-\alpha)},
%$$
$$
\theta=\min_{x,y\in [m,M], \alpha\in[\frac{1}{2},1]}\frac{\log\bigg(\frac{\log(\alpha x+(1-\alpha)y)}{\log(xy)}\bigg)}{\log(1-\alpha)}
$$
and choose $s>0$ such that $s\leq \theta$ to deduce the result.
\end{proof}

\begin{corollary}\label{cor-1}
Let $A$ be a positive operator on a Hilbert space
$\mathcal H$ satisfying $1< A\leq M$.
Then the following inequality holds for every unit vector $x$, $q\geq\frac{M}{M-1}$, and for some $0<s<1$: %, every $q\geq\frac{e}{e-1}$,
\begin{equation}\label{eq-rem0}
\langle \log(A)x,x\rangle \leq  \frac{2^{1-s}\log(M)}{(q-1)(M-1)}\bigg(\frac{q-1}{q}\bigg)^q\langle Ax,x\rangle^{q}.
\end{equation}
%where $q>1$, $\frac{\log(m)}{m}q\leq\frac{\log(M)-\log(m)}{M-m}\leq \frac{\log(M)}{M}q$.
%$f(M)>f(m)$, $\frac{f(M)}{M}<\frac{f(m)}{m}$ and the condition $\frac{f(M)}{M}q\leq\frac{f(M)-f(m)}{M-m}\leq\frac{f(m)}{m}q$
%holds for a real number $q$ such that $0<s<q<1$.
\end{corollary}
\begin{proof}
Put $f(t)=\log(t)$ for $1< t\leq M$ in Theorem \ref{main-0}(i).
By Lemma \ref{lemma-log} the function $f(t)$ is a real valued $s$-convex function for some $0<s<1$
on $(1, M]$.
The function $f$ satisfies the conditions of Theorem \ref{main-0}(i) on $(1,M]$.
Indeed, since $f$ is increasing, $\log(M)>\log(1)=0$.
So, $\frac{\log(M)}{M}>0$ and
by assumption since $q\geq \frac{M}{M-1}$,
one has
$$
0\leq \frac{\log(M)}{M-1}\leq \frac{\log(M)}{M}q.
$$
According to Theorem \ref{main-0}(i), we find that
\begin{equation}\label{eq-rem0}
\langle \log(A)x,x\rangle \leq 2^{1-s} K_{\log}(1,M,q) \langle Ax,x\rangle^{q}.
\end{equation}
Note that $K_{\log}(1,M,q)=\frac{\log(M)}{(q-1)(M-1)}(\frac{q-1}{q})^q$.
%Since $A\geq 1$, by taking minimum over $q\geq \frac{e}{e-1}$, we reach the result.
Whence the proof is complete by Theorem \ref{main-0}(i).
\end{proof}

\begin{remark}
We remark that Corollary \ref{cor-1} is not a trivial result.
The following inequality
is due to Mond--Pe$\check{c}$ari$\acute{c}$ \cite{Mond-HJM1993}
\begin{equation}\label{eq-rem1}
\langle \log(A)x,x\rangle\leq \log(\langle Ax,x\rangle)
\end{equation}
for every positive operator $A$ and every unit vector $x$.
%On the other hand, in Corollary
By a numerical computation we compare two upper bounds provided in \eqref{eq-rem0} and \eqref{eq-rem1}.
We show that there is no ordering between the upper bounds given in the inequalities \eqref{eq-rem0} and \eqref{eq-rem1}
and so Corollary \ref{cor-1} is not a trivial result.

Consider
$A=\Big(
\begin{array}{ll}
1.0 & 0.0\\
0.0 & 1.1
\end{array}
 \Big)$,
  $M=10$, $q=2$, and
 $x=\Bigg(
\begin{array}{l}
\frac{1}{\sqrt{2}}\\
\frac{1}{\sqrt{2}}
\end{array}
 \Bigg)$.
 Let $\log t$ be the decimal or common logarithm.
So, one can see
\begin{align*}
\langle Ax,x\rangle&=1.05,\\
\log(\langle Ax,x\rangle)&=0.021189,\\
2^{1-s} K_{\log}(1,M,q)\langle Ax,x\rangle^{q}&=2^{1-s}(0.030625) > 0.021189.
\end{align*}
This entails
\begin{equation}\label{eq-comp}
2^{1-s} K_{\log}(1,M,q)\langle Ax,x\rangle^{q} > \log(\langle Ax,x\rangle),
\end{equation}
while for $q=8$ a numerical computation shows
$$
2^{1-s} K_{\log}(1,M,q)\langle Ax,x\rangle^{q}<(2)(0.00806) < 0.021189
$$
and so the reverse inequality holds in \eqref{eq-comp}.
\end{remark}

%%%%%%%%%%%%%%%%%%%%%%%%%%%%%%%%%%%%%%%%%%%%%%%%%%%%%%%%%%%%%%%%%%%%%%%%%%%%%%%%
We provide a complementary result for Theorem \ref{main-0} in the case $0<q<1$.

\begin{theorem}\label{main-1}
Let $A$ be a positive operator on a Hilbert space
$\mathcal H$ satisfying $M\geq A\geq m>0$.
Also let $f(t)$ be a real valued continuous $s$-convex function on $[m,M]$.
Then the following inequality holds for every unit vector $x$:
$$
\langle f(A)x,x\rangle\leq 2^{1-s}\max\Big\{\frac{f(m)}{m^q},\frac{f(M)}{M^q}\Big\}\langle Ax,x\rangle^{q},
$$
where
$f(M)>f(m)$, $\frac{f(M)}{M}<\frac{f(m)}{m}$ and the condition $\frac{f(M)}{M}q\leq\frac{f(M)-f(m)}{M-m}\leq\frac{f(m)}{m}q$
holds for a real number $q$ such that $0<q<1$.
\end{theorem}
\begin{proof}
As in the proof of Theorem \ref{main-0}, one has
\begin{align}\label{eq-qq500}
\langle f(A)x,x\rangle\leq 2^{1-s}\Big(\max_{m\leq t\leq M}h(t)\Big)\langle Ax,x\rangle^{q}.
\end{align}
Putting $K=f(M)$ and $k=f(m)$ in Lemma \ref{lemma2}, then the conditions in Theorem \ref{main-1}
corresponds to the conditions in Lemma \ref{lemma2}, so the proof is complete by \eqref{eq-qq500}.
\end{proof}

\begin{corollary}\label{cor-2}
Let $A$ be a positive operator on a Hilbert space
$\mathcal H$ satisfying $M\geq A\geq m>0$.
Then the following inequality holds for every unit vector $x$:
\begin{equation}\label{eq-rem10}
\langle A^px,x\rangle\leq 2^{1-s}\max\{m^{p-q},M^{p-q}\}\langle Ax,x\rangle^{q},
\end{equation}
where
$0<s<p<1$, $0<q<1$, and $M^{p-1}q\leq\frac{M^{p}-m^{p}}{M-m}\leq m^{p-1}q$.
\end{corollary}
\begin{proof}
Put $f(t)=t^p$ for $0<s<p<1$ in Theorem \ref{main-1}.
By Lemma \ref{lm-scnx} the function $f(t)$ is a real valued continuous $s$-convex function
on $[m, M]$. On the other hand, since $f$ is increasing, $M^p>m^p$.
For $g(t)=t^{p-1}$, $0<s<p<1$, the function $g(t)$
is decreasing on $[m, M]$ and so $M^{p-1}<m^{p-1}$
holds for $0<s<p<1$, that is, $f(M)=M^p>m^p=f(m)$ and $\frac{f(M)}{M}=M^{p-1}<m^{p-1}=\frac{f(m)}{m}$.
Whence the proof is complete by Theorem \ref{main-1}.
\end{proof}

%%%%%%%%%%%%%%%%%%%%%%%%%%%%%%%%%%%%%%%%%%%%%%%%%%%%%%%%%%%%%%%%%%%%%%%%%
%remark
%%%%%%%%%%%%%%%%%%%%%%%%%%%%%%%%%%%%%%%%%%%%%%%%%%%%%%%%%%%%%%%%%%%%%%%%%%%%%%%%

\begin{corollary}\label{cor-2020}
Let $A$ be a positive operator on a Hilbert space
$\mathcal H$ satisfying $M\geq A\geq m>0$.
Then the following inequality holds for every unit vector $x$:
\begin{equation}%\label{eq-rem10}
||A||^p\leq 2^{1-s}\max\{m^{p-q},M^{p-q}\}\langle Ax,x\rangle^{q},
\end{equation}
where
$0<s<p<1$, $0<q<1$, and $M^{p-1}q\leq\frac{M^{p}-m^{p}}{M-m}\leq m^{p-1}q$.
\end{corollary}
\begin{proof}
Consider $f(t)=||t||^p$ for $0<s<p<1$ in Theorem \ref{main-1}.
The function $f(t)$ is a real valued continuous $s$-convex function
on $[m, M]$ by Lemma \ref{lm-scnx-norm}.
The rest of the proof is similar to that of Corollary \ref{cor-2}.
\end{proof}

\begin{corollary}\label{cor-2021}
Let $A$ be a positive operator on a Hilbert space
$\mathcal H$ satisfying $M\geq A\geq m>0$.
Then the following inequality holds for every unit vector $x$:
\begin{equation}%\label{eq-rem10}
\frac{1}{2^{\frac{1-s}{p}}}||A||\leq \langle Ax,x\rangle \leq ||A||,
\end{equation}
where
$0<s<p<1$ and $M^{p-1}p\leq\frac{M^{p}-m^{p}}{M-m}\leq m^{p-1}p$.
\end{corollary}
\begin{proof}
The first inequality comes from Corollary \ref{cor-2020} by letting $q=p$.
The second inequality follows from the Cauchy-Schwartz inequality, i.e.,
$$
\langle Ax,x\rangle\leq ||Ax||||x||\leq ||A||||x||^2=||A||
$$
for every unit vector $x$ and every positive operator $A$.
\end{proof}

%%%
Let $A$ be a bounded linear operator on a Hilbert space $\mathcal H$.
The numerical range $W(A)$ is defined by
$$
W(A) := \{\langle Ax, x\rangle : x \in \mathcal H, ||x|| = 1\}.
$$
The largest absolute value of the numbers in the numerical range is called the numerical radius, i.e.
$$
w(A) := \sup\{|\lambda|: \lambda\in W(A)\}=\sup_{||x|| = 1}|\langle Ax, x\rangle|.
$$

It is known that $W(A)$ lies in the closed disc of radius $||A||$ centered at the origin
and so $w(A) \leq ||A||$.
The following corollary is a simple consequence of Corollary \ref{cor-2020}.

\begin{corollary}\label{cor-202100}
Let $A$ be a positive operator on a Hilbert space
$\mathcal H$ satisfying $M\geq A\geq m>0$.
Then the numerical range of $A$ is contained in the closed interval $\bigg[\frac{1}{2^{\frac{1-s}{p}}}||A||,||A||\bigg]$
and
$$
\frac{1}{2^{\frac{1-s}{p}}}||A||\leq w(A) \leq ||A||,
$$
where
$0<s<p<1$, $1<p+s$, and $M^{p-1}p\leq\frac{M^{p}-m^{p}}{M-m}\leq m^{p-1}p$.
\end{corollary}

%%%%%%%%%%%%%%%%%%%%%%%%%%%%%%%%%%%%%%%%%%%%%%%%%%%%%%%%%%%%%%
The advantages of this corollary is that it gives a better lower bound for the numerical radius of
the positive operator $A$
under some conditions than the lower bound presented in \cite[Theorem 1.2]{Goldberg-Tadmor1982}.
Indeed, it has been proved in \cite[Theorem 1.2]{Goldberg-Tadmor1982} that
$\frac{1}{2}||A||\leq w(A)$, whence under the conditions of this corollary one has
$$
\frac{1}{2}||A||<2^{\frac{s-1}{p}}||A||\leq w(A).
$$

We refined in \cite[Corollary 5]{Nikoufar-CAOT2021} the H\"{o}lder-McCarthy inequality by providing an upper bound as follows.
We apply this corollary in our results.

\begin{corollary}\label{Mac1}
Let $A$ be a positive operator on a Hilbert space $\mathcal H$. Then,
\begin{itemize}
\item[(i)] for all $0<s\leq q<1$ and $x\in \mathcal H$ with $||x||=1$,
\begin{align}\label{eq-corMac1}
\langle A^qx,x\rangle\leq \langle Ax,x\rangle^q\leq 2^{1-s}\langle A^qx,x\rangle,
\end{align}
\item[(ii)] for all $\frac{1}{s}\geq q>1$ and $x\in \mathcal H$ with $||x||=1$,
\begin{align}\label{eq-corMac2}
\langle Ax,x\rangle^q\leq \langle A^qx,x\rangle\leq 2^{(1-s)q}\langle Ax,x\rangle^q.
\end{align}
\end{itemize}
\end{corollary}

As another application of our results, we provide
a complementary result of an order preserving inequality of three variables
given by Mi$\acute{c}$i$\acute{c}$--Pe$\check{c}$ari$\acute{c}$--Seo \cite{Micic-LAA2003}.
The following result associated with the H\"{o}lder-McCarthy inequality.

\begin{corollary}\label{cor-3}
Let $A$ and $B$ be positive operators on a Hilbert space
$\mathcal H$ satisfying $m\leq A\leq M$ with $0<m<M$ and $A\leq B$.
Then %the following inequality holds for every unit vector $x$:
\begin{equation}%\label{eq-rem11}
A^p\leq 2^{2(1-s)}\max\{m^{p-q},M^{p-q}\} B^{q},
\end{equation}
where
$0<s<p<1$, $0<s\leq q<1$, and $M^{p-1}q\leq\frac{M^{p}-m^{p}}{M-m}\leq m^{p-1}q$.
\end{corollary}
\begin{proof}
By using Corollary \ref{cor-2} and applying the refined H\"{o}lder-McCarthy inequality \eqref{eq-corMac1} one can see that
\begin{align*}
\langle A^px,x\rangle&\leq 2^{1-s}\max\{m^{p-q},M^{p-q}\}\langle Ax,x\rangle^{q}\\
&\leq 2^{1-s}\max\{m^{p-q},M^{p-q}\}\langle Bx,x\rangle^{q}\\
&\leq 2^{2(1-s)}\max\{m^{p-q},M^{p-q}\}\langle B^{q}x,x\rangle
\end{align*}
for every unit vector $x$ and so the result follows.
\end{proof}

\section{A difference version of the Kantorovich type inequalities}
Recently, Furuta and Giga \cite{Furuta-Alg2003} gave the complementary result of the Kantorovich
type order preserving inequalities by Mi$\acute{c}$i$\acute{c}$--Pe$\check{c}$ari$\acute{c}$--Seo.
In this section, we show a difference version of Furuta and Giga's result for $s$-convex functions.
We first prove the following lemma which is a difference version of Lemma \ref{lemma2}.
%%%%%%%%%%%%%%%%%%%%%%%%%%%%%%%%%%%%%%%%%%%%%%%%%%%%%%%%%%%%%%%%%%%%%%%%%%%%%%%%
\begin{lemma}\label{lemma1-difference}
Let $u(t)$ be defined by
\begin{equation}\label{eq-u}
u(t)=k+\frac{K-k}{M-m}(t-m)-t^q
\end{equation}
on $[m,M]$, $(M>m>0)$,
and any real number $q$ such that $0<q<1$
and any real numbers $K$ and $k$. Then
the function $u(t)$ has the following lower bound
on $[m,M]$:
\begin{equation}\label{eq-ppp}
k+\frac{K-k}{M-m}\bigg(\bigg(\frac{K-k}{q(M-m)}\bigg)^{\frac{1}{q-1}}-m\bigg)-\bigg(\frac{K-k}{q(M-m)}\bigg)^{\frac{q}{q-1}}
\end{equation}
where $K>k$ and $qm^{q-1}\geq\frac{K-k}{M-m}\geq qM^{q-1}$.
Moreover, in this case the upper bound occurs at the beginning or end of the interval $[m,M]$
and the upper bound value is $\max\{k-m^q, K-M^q\}$.
\end{lemma}
\begin{proof}
By an easy differential calculus, we know that $u'(t_1)=0$
when
$$
t_1=\bigg(\frac{K-k}{q(M-m)}\bigg)^{\frac{1}{q-1}}.
$$
On the other hand, $u''(t_1)>0$ and so $t_1$ gives
the lower bound \eqref{eq-ppp} of $u(t)$ on $[m,M]$.
Moreover, $t_1$ is the global minimum point of $u$
and $m$ and $M$ are the local maximum point of $u$ with the maximum values $k-m^q$ and $K-M^q$, respectively.
\end{proof}

%%%%%%%%%%%%%%%%%%%%%%%%%%%%%%%%%%%%%%%%%%%%%%%%%%%%%%%%%%%%%%%%%%%%%%%%%%%%%%%%
\begin{lemma}\label{lemma2-difference}
Let $u(t)$ be defined by \eqref{eq-u} on $[m,M]$, $(M>m>0)$,
and any real number $q$ and any real numbers $K$ and $k$.
Then the function $u(t)$ has the upper bound \eqref{eq-ppp}
on $[m,M]$, where $m, M, k, K$ and $q$
satisfy any one (i) or (ii) respectively:
\begin{itemize}
\item[(i)] $K>k$ and $qm^{q-1}\leq\frac{K-k}{M-m}\leq qM^{q-1}$ holds for a real number $q>1$,
\item[(ii)] $K<k$ and $qM^{q-1}\leq\frac{K-k}{M-m}\leq qm^{q-1}$ holds for a real number $q<0$.
\end{itemize}
\end{lemma}
\begin{proof}
The results follow from the fact that in these cases
$u''(t_1)<0$.
\end{proof}

%It is remarkable that Theorem \ref{CAOT-thm1} holds for $q<0$ or $q>1$.
%We now extend Theorem \ref{CAOT-thm1} for the case $0<q<1$ as a complementary result.

\begin{theorem}\label{thm1-diff}
Let $A$ be a positive operator on a Hilbert space
$\mathcal H$ satisfying $M\geq A\geq m>0$.
Also let $f(t)$ be a real valued continuous $s$-convex function on $[m,M]$.
Then the following inequality holds for every unit vector $x$:
$$
\langle f(A)x,x\rangle\leq 2^{1-s}\bigg(\max\{f(m)-m^q,f(M)-M^q\}+\langle Ax,x\rangle^{q}\bigg),
$$
where
$f(M)>f(m)$ and the condition $qm^{q-1}\geq\frac{K-k}{M-m}\geq qM^{q-1}$
holds for a real number $q$ such that $0<q<1$.
\end{theorem}
\begin{proof}
Subtracting $2^{1-s}\langle Ax,x\rangle^{q}$ from both sides of \eqref{eq-qq2}, we get
\begin{align}\label{eq1-diff}
\langle f(A)x,x\rangle-2^{1-s}\langle Ax,x\rangle^{q}
\leq 2^{1-s} u(\langle Ax,x\rangle),
\end{align}
where
\begin{align}\label{eq2-diff}
u(t)=f(m)+\frac{f(M)-f(m)}{M-m}(t-m)-t^{q}
\end{align}
for every $t\in [m,M]$.
So, one can deduce from \eqref{eq1-diff} that
\begin{align}\label{eq3-diff}
\langle f(A)x,x\rangle-2^{1-s}\langle Ax,x\rangle^{q}
\leq 2^{1-s}\Big(\max_{m\leq t\leq M}u(t)\Big).
\end{align}
Putting $K=f(M)$ and $k=f(m)$ in Lemma \ref{lemma1-difference}, then the conditions in Theorem \ref{thm1-diff}
corresponds to the conditions in Lemma \ref{lemma1-difference}, so the result follows from \eqref{eq3-diff}.
\end{proof}

\begin{theorem}\label{thm2-diff}
Let $A$ be a positive operator on a Hilbert space
$\mathcal H$ satisfying $M\geq A\geq m>0$.
Also let $f(t)$ be a real valued continuous $s$-convex function on $[m,M]$.
Then the following inequality holds for every unit vector $x$ and for any
real number $q$ depending on (i) or (ii):
$$
\langle f(A)x,x\rangle\leq 2^{1-s} (K_f^d(m,M,q)+\langle Ax,x\rangle^q),
$$
where
\begin{align*}
K_f^d(m,M,q)=f(m)&+\frac{f(M)-f(m)}{M-m}\bigg(\bigg(\frac{f(M)-f(m)}{q(M-m)}\bigg)^{\frac{1}{q-1}}-m\bigg)\\
&-\bigg(\frac{f(M)-f(m)}{q(M-m)}\bigg)^{\frac{q}{q-1}}
\end{align*}
under any one of the following conditions (i) or (ii) respectively:
\begin{enumerate}
\item[(i)] $f(M)>f(m)$ and $qm^{q-1}\leq\frac{f(M)-f(m)}{M-m}\leq qM^{q-1}$
holds for a real number $q>1$,
\item[(ii)] $f(M)<f(m)$ and $qM^{q-1}\leq\frac{f(M)-f(m)}{M-m}\leq qm^{q-1}$
holds for a real number $q<0$.
\end{enumerate}
\end{theorem}
\begin{proof}
By the same approach as in the proof of Theorem \ref{thm1-diff}
and using Lemma \ref{lemma2-difference} one can reach the result.
\end{proof}

\begin{corollary}\label{cor1-diff}%cor-2
Let $A$ be a positive operator on a Hilbert space
$\mathcal H$ satisfying $M\geq A\geq m>0$.
Then the following inequality holds for every unit vector $x$:
\begin{equation}%\label{eq-rem10}
\langle A^px,x\rangle\leq 2^{1-s}\bigg(\max\{m^{p-q},M^{p-q}\}+\langle Ax,x\rangle^{q}\bigg),
\end{equation}
where
$0<s<p<1$, $0<q<1$, and $M^{q-1}q\leq\frac{M^{p}-m^{p}}{M-m}\leq m^{q-1}q$.
\end{corollary}
\begin{proof}
Put $f(t)=t^p$ for $0<s<p<1$ in Theorem \ref{thm1-diff}.
Like what was mentioned in Corollary \ref{cor-2}, the function $f$ satisfies the conditions of Theorem \ref{thm1-diff}.
Whence the proof is complete by Theorem \ref{thm1-diff}.
\end{proof}

As an application of this result, we give
a complementary result of an order preserving inequality of three variables.

\begin{corollary}%\label{cor-3}
Let $A$ and $B$ be positive operators on a Hilbert space
$\mathcal H$ satisfying $m\leq A\leq M$ with $0<m<M$ and $A\leq B$.
Then %the following inequality holds for every unit vector $x$:
\begin{equation}%\label{eq-rem11}
A^p\leq 2^{1-s}\bigg(\max\{m^{p-q},M^{p-q}\} + 2^{1-s} B^{q}\bigg),
\end{equation}
where
$0<s\leq p<1$, $0<s\leq q<1$, and $M^{q-1}q\leq\frac{M^{p}-m^{p}}{M-m}\leq m^{q-1}q$.
\end{corollary}
\begin{proof}
By using Corollary \ref{cor1-diff} and the refined H\"{o}lder-McCarthy inequality \eqref{eq-corMac1} one can get the result,
like what was done in Corollary \ref{cor-3}.
\end{proof}

%{\bf Acknowledgement:}
%The first author was supported by Payame Noor University under grant number~46384.

\section{Declarations}
We remark that the potential conflicts of interest and
data sharing not applicable to this article and no data sets were generated during the current study.

\begin{acknowledgements}
We would like to express our sincere thanks for the referee for the useful comments
that improved the manuscript.
The first author was supported by Payame Noor University.% under grant number~46384.
\end{acknowledgements}


\begin{thebibliography}{99}
%----------------------------------------------------------------------%
%\bibitem{Agarwal-Dragomir-CMA2010}
%R. P. Agarwal and S. S. Dragomir,
%A survey of Jensen type inequalities for functions of self-adjoint operators in Hilbert spaces,
%Comput. Math. Appl. 59(12), 3785--3812 (2010).
%----------------------------------------------------------------------%
%\bibitem{Alirezaei2018}
%G. Alirezaei and R. Mazhar,
%On exponentially concave functions and their impact in information theory,
%J. Inform. Theory Appl. 9(5), 265--274 (2018).
%----------------------------------------------------------------------%
%\bibitem{Bernstein1929}
%S. N. Bernstein,
%Sur les fonctions absolument monotones,
%Acta Math. 52, 1--66 (1929).
%----------------------------------------------------------------------%
%\bibitem{Dragomir-SJM1995}
%S. S. Dragomir, J. Pe$\check{c}$ari$\acute{c}$, and L. E. Persson,
%Some inequalities of Hadamard type,
%Soochow J. Math. 21(3), 335--341 (1995).
%----------------------------------------------------------------------%
%\bibitem{Dragomir-BAMS1998}
%S. S. Dragomir and C. E. M. Pearce,
%Quasi-convex functions and Hadamard's inequality,
%Bull. Aust. Math. Soc. 57(3), 377--385 (1998).
%----------------------------------------------------------------------%
\bibitem{Cvetkovski}
Z. Cvetkovski,
Inequalities: Theorems, Techniques and Selected Problems;
Springer: Berlin/Heidelberg, Germany, 2012; ISBN 978-3-642-23791-1.
%----------------------------------------------------------------------%
\bibitem{Dragomir-DM1999}
S. S. Dragomir and S. Fitzpatrick,
The Hadamard's inequality for $s$-convex functions in the second sense,
Demonstratio Math. 32(4), 687--696 (1999).
%----------------------------------------------------------------------%
\bibitem{Dragomir-book2002}
S. S. Dragomir and Th. M. Rassias,
Ostrowski type inequalities and applications in numerical integration,
Kluwer Academic Publishers, Dordrecht, Boston, London, 2002.
%----------------------------------------------------------------------%
\bibitem{Kirmaci-AMC2007}
U. S. Kirmaci, M. Klari$\check{c}$i$\acute{c}$, M. E. $\ddot{O}$zdemirc, and J. Pe$\check{c}$ari$\acute{c}$,
Hadamard--type inequalities for $s$-convex functions,
Appl. Math. Comput. 193, 26--35 (2007).
%----------------------------------------------------------------------%
\bibitem{Furuta-JIA1998}
T. Furuta,
Operator inequality associated with H$\ddot{o}$lder-McCarthy and Kantorovich inequalities,
J. of lnequal. and Appl. 2, 137--148 (1998).
%----------------------------------------------------------------------%
\bibitem{Furuta-TF2001}
T. Furuta,
Invitation to Linear Operators,
Taylor and Francis, London, 2001.
%----------------------------------------------------------------------%
\bibitem{Furuta-Alg2003}
T. Furuta and M. Giga,
A complementary result of Kantorovich type order preserving inequalities
by Mi$\acute{c}$i$\acute{c}$--Pe$\check{c}$ari$\acute{c}$--Seo,
Linear Alg. Appl. 369, 27--40 (2003).
%----------------------------------------------------------------------%
\bibitem{Furuta-Zagreb2005}
T. Furuta, J. Mi$\acute{c}$i$\acute{c}$, J. Pe$\check{c}$ari$\acute{c}$, and Y. Seo,
Mond--Pe$\check{c}$ari$\acute{c}$ method in operator inequalities,
inequalities for bounded self-adjoint operators on a Hilbert Space,
Element, Zagreb, 2005.
%----------------------------------------------------------------------%
\bibitem{Goldberg-Tadmor1982}
M. Goldberg and E. Tadmor,
On the numerical radius and its applications,
Linear Alg. Appl. 42, 263--284 (1982).
%----------------------------------------------------------------------%
\bibitem{Hudzik-AEqMath1994}
H. Hudzik and L. Maligranda,
Some remarks on $s$-convex functions,
Aequationes Math. 17, 100--111 (1994).
%----------------------------------------------------------------------%
\bibitem{Mond-HJM1993}
B. Mond and J. Pe$\check{c}$ari$\acute{c}$,
Convex inequalities in Hilbert space,
Houston J. Math. 19, 405--420 (1993).
%----------------------------------------------------------------------%
\bibitem{Mond-IJM1993}
B. Mond, J.E. Pe$\check{c}$ari$\acute{c}$,
Convex inequalities for several self--adjoint operators on a Hilbert space,
Indian J. Math. 35, 121--135 (1993).
%----------------------------------------------------------------------%
\bibitem{Micic-LAA2003}
J. Mi$\acute{c}$i$\acute{c}$, J. Pe$\check{c}$ari$\acute{c}$, and Y. Seo,
Function order of positive operators based on the Mond--Pe$\check{c}$ari$\acute{c}$ method,
Linear Alg. Appl. 360, 15--34 (2003).
%----------------------------------------------------------------------%
\bibitem{Nikoufar-Filomat2020}
I. Nikoufar and D. Saeedi,
Some Inequalities for $P$-Class Functions,
Filomat 34(13), 4555--4566 (2020).
%----------------------------------------------------------------------%
\bibitem{Nikoufar-CAOT2021}
I. Nikoufar and D. Saeedi,
An operator version of the Jensen inequality for $s$-convex functions,
Complex Anal. Oper. Theory 15, 92 (2021). https://doi.org/10.1007/s11785-021-01139-x
%----------------------------------------------------------------------%
%\bibitem{Noor2020}
%M. A. Noor,  K. I. Noor, and M. Th. Rassias,
%New Trends in General Variational Inequalities,
%Acta Applicandae Mathematicae, 170(1), 981--1064 (2020).
%----------------------------------------------------------------------%
%\bibitem{Pearce-Dragomir}
%C. E. M. Pearce and S. S. Dragomir,
%Selected topics on Hermite-Hadamard inequalities and applications
%[https://rgmia.org/papers/monographs/Master.pdf].
%----------------------------------------------------------------------%
\bibitem{Pinheiro-AEqMath2007}
M. R. Pinheiro,
Exploring the concept of $s$-convexity,
Aequationes Math. 74, 201--209 (2007).
%----------------------------------------------------------------------%
\bibitem{Set-CA2012}
E. Set,
New inequalities of Ostrowski type for mappings whose derivatives
are $s$-convex in the second sense via fractional integrals,
Appl. Math. Comput. 63, 1147--1154 (2012).
%----------------------------------------------------------------------%
\bibitem{Takahashi-MJ1999}
S.-E. Takahashi, M. Tsukada, K. Tanahashi, T. Ogiwara,
An inverse type of Jensen's inequality,
Math. Japon. 50, 85--92, (1999).
\end{thebibliography}


\end{document}